\documentclass[12pt, reqno]{amsart}
\usepackage{amssymb, amsmath, amsthm, cancel, hyperref, color, enumerate, verbatim, mathtools}

\usepackage[margin=1in, includeheadfoot]{geometry}

\newtheorem{thm}{Theorem}
\newtheorem{rem}{Remark}
\newtheorem{ex}{Example}

\renewcommand{\Re}{\mathrm{Re}}
\newcommand{\SL}{\mathrm{SL}}
\newcommand{\PGL}{\mathrm{PGL}}

\begin{document}

\title[Zeta-polynomials]{Zeta-polynomials, Hilbert polynomials, and the Eichler-Shimura identities}
\author{Marie Jameson}

\maketitle
\section*{Abstract}
In 2017, Ono, Rolen, and Sprung \cite{OnoRolenSprung17} answered problems of Manin \cite{Manin16} by defining zeta-polynomials $Z_f(s)$ for even weight newforms $f\in S_k(\Gamma_0(N)$; these polynomials can be defined by applying the ``Rodriguez-Villegas transform'' to the period polynomial of $f$. It is known that these zeta-polynomials satisfy a functional equation $Z_f(s) = \pm Z_f(1-s)$ and they have a conjectural arithmetic-geometric interpretation. Here, we give analogous results for a slightly larger class of polynomials which are also defined using the Rodriguez-Villegas transform.

\section{Introduction and statement of results}\label{intro}

Let $f\in S_k(\Gamma_0(N)$ be a newform of even weight $k$ and level $N$, and let $L(f, s)$ be the $L$-function associated to $f$. Manin \cite{Manin16} speculated that the critical $L$-values
\[L(f,1), L(f,2), \ldots, L(f,k-1)\]
can be assembled in a natural way to build a \emph{zeta-polynomial}. This polynomial $Z_f(s)$ should
\begin{enumerate}[(i)]
\item satisfy a functional-equation $Z_f(s) = \pm Z_f(1-s)$,
\item obey the ``Riemann hypothesis:'' if $Z_f(\rho)=0$ then $\Re(\rho)=1/2,$ and 
\item have an arithmetic-geometric interpretation.
\end{enumerate}
Recently, Ono, Rolen, and Sprung \cite{OnoRolenSprung17} defined a zeta-polynomial $Z_f(s)$ which satisfies properties (i) and (ii) above.  Assuming the Bloch–Kato Tamagawa Number Conjecture, it also satisfies property (iii) by encoding the arithmetic of a combinatorial arithmetic-geometric object called “Bloch–Kato complex” for $f$.

Although $Z_f(s)$ can be defined as a sum involving weighted moments of critical $L$-values and signed Stirling numbers of the first kind, it is more convenient here to instead express it in terms of (a slightly normalized version of) the \emph{period polynomial} of $f$, which is given by
\[R_f(X) := 
(\sqrt{N})^{k-1}\frac{(k-2)!}{(2\pi)^{k-1}}\sum_{n=0}^{k-2}\frac{(2\pi X)^n}{n!(\sqrt{N})^n}L(f,k-n-1).\]
Period polynomials are well-studied objects which are known to have many beautiful properties.  For example, it is known that $R_f(X)$ satisfies its own Riemann hypothesis: all of its roots occur on the unit circle $|X|=1$, as proved in \cite{ElGuindyRaji14, JinMaOnoSound16}.  See Section \ref{perpolybackground} for additional background information about $R_f(X).$ In \cite{OnoRolenSprung17}, $Z_f(s)$ is described as the unique polynomial which satisfies
\[\frac{R_f(X)}{(1-X)^{w+1}} = \sum_{n\geq 0}Z_f(-n)X^n,\]
where $w:=k-2.$ This relationship between $R_f(X)$ and $Z_f(s)$ is known as the ``Rodriguez-Villegas transform,'' and a key theorem of Rodriguez-Villegas \cite{RodriguezVillegas02} allows Ono, Rolen, and Sprung to translate the Riemann hypothesis for $R_f(X)$ into statements (i) and (ii) about $Z_f(s).$

Zeta-polynomials are relatively new objects, and little else is currently known about their properties.  However, the results described thus far give evidence that known properties of a newform $f$ and its period polynomial $R_f(X)$ could be translated into the realm of the zeta-polynomial $Z_f(s).$ This could give us more insight into the behavior of zeta-polynomials.  The goal of this article is to offer additional evidence in this direction.

Note, however, that the results here are general enough that they do not require $R(X)$ to be the period polynomial of a newform (and for example, our results will apply to even/odd parts of period polynomials).  Thus we fix the following notation: let $w\geq 2$ be even, let $R(X) \in \mathbb{C}[X]$ be \emph{any} polynomial of degree at most $w$, and let $Z(s)$ be the unique polynomial which satisfies
\begin{equation}\label{defofZ}
\frac{R(X)}{(1-X)^{w+1}} = \sum_{n\geq 0}Z(-n)X^n,
\end{equation}

Our first result assumes the identity
\begin{equation}\label{fricke}
R(X) + \varepsilon i^{w}X^wR(1/X) = 0
\end{equation}
(where $\varepsilon=\pm 1$ is a constant) and interprets its meaning in terms of $Z(s).$ Note here that equation \eqref{fricke} is important because is it known to be true for any period polynomial (as well as its even and odd parts) associated to a newform $f\in S_k(\Gamma_0(N)),$ where $\varepsilon = \pm 1$ is the eigenvalue of $f$ under the Fricke involution.

\begin{thm}\label{eichshim1}
Let $w\geq 2$ be even, let $R(X) \in \mathbb{C}[X]$ be any polynomial of degree at most $w$, and let $Z(s)$ be the polynomial satisfying \eqref{defofZ}. If $R(X)$ satisfies \eqref{fricke}, then we have that
\[Z(s) +\varepsilon i^wZ(1-s) =0.\]
\end{thm}

\begin{rem}
Since the period polynomial $R_f(X)$ satisfies equation \eqref{fricke} and the conclusion of the theorem gives property (i) above, one can view Theorem \ref{eichshim1} as a generalization of property (i) proved by Ono, Rolen, and Sprung \cite{OnoRolenSprung17} which uses completely different techniques (and does not depend on the Riemann hypothesis for period polynomials).
\end{rem}

One may also consider what the ``Eichler-Shimura relations'' for period polynomials associated to cusp forms of level 1 (as described in Section \ref{prelim}) tell us about zeta-polynomials.  Thus we suppose that
\begin{align} 
R(X) + (-iX)^w R(1/X)&= 0 \label{Res1}\\
R(X) + (-iX)^wR\left(\frac{X-i}{-iX}\right) + (-iX-1)^wR\left(\frac{-i}{-iX-1}\right) &= 0, \label{Res2}
\end{align}
and obtain the following result.

\begin{thm}\label{eichshim2}
Let $w\geq 2$ be even, let $R(X) \in \mathbb{C}[X]$ be any polynomial of degree at most $w$, and let $Z(s)$ be the polynomial satisfying \eqref{defofZ}. If $R(X)$ satisfies \eqref{Res1} and \eqref{Res2}, then
\[Z(s) + i^wZ(1-s) =0\]
and for any positive integer $n$ we have that
\begin{equation*} \begin{split}
Z(-n) + &(-i)^w\sum_{m=1}^{n+1} a_{-m}Z(1-m) +\\
&\sum_{k\geq 0}\sum_{m=0}^{k+n} \sum_{j=0}^{k+n-m} \binom{k+n}{n} \binom{m+w}{w} \binom{w+1}{j} \frac{(-1)^{j+1}(-i)^k }{(1-i)^{m+w+1}}  Z(m+j-k-n) = 0
\end{split}\end{equation*}
where
\[\frac{(1-x)^{w+1} (x+i)^{n}}{(x+i-ix)^{w+1}(ix)^{n+1}} = \sum_{m=-n-1}^\infty a_mx^m.\]
\end{thm}

\begin{ex}
For instance, consider the unique newform $\Delta \in S_{12}(\Gamma_0(1)).$ As computed in \cite{OnoRolenSprung17}, we have that
\begin{equation*} \begin{split}
R_\Delta(X) \approx 0.114379\cdot\left(\frac{36}{691}X^{10} + X^8 + 3X^6 + 3X^4 + X^2 + \frac{36}{691}\right)\\
+ 0.00926927\cdot\left(4X^9 + 25X^7 + 42X^5 + 25X^3 + 4X\right)
\end{split} \end{equation*}
and
\begin{equation*} \begin{split}
Z_\Delta(s) \approx (5.11\times 10^{-7})s^{10} - (2.554\times 10^{-6})s^9 + (6.01\times 10^{-5})s^8 - (2.25\times 10^{-4})s^7\\
 + 0.00180s^6 - 0.00463s^5 + 0.0155s^4 - 0.0235s^3 + 0.0310s^2 - 0.0199s + 0.00596.
\end{split} \end{equation*}
Since $\Delta$ is the unique normalized cusp form of weight 12 and level 1, we have that $R_\Delta(X)$ satisfies equations \eqref{fricke}, \eqref{Res1}, and \eqref{Res2}, so Theorems \ref{eichshim1} and \ref{eichshim2} apply, although of course the first statement of Theorem \ref{eichshim2} that
\[Z_\Delta(s) - Z_\Delta(1-s) =0\]
is the same as Theorem \ref{eichshim1}, and it is already known by \cite{OnoRolenSprung17}.

However, one can now consider
\[R_\Delta^-(X) := 4X^9 + 25X^7 + 42X^5 + 25X^3 + 4X\]
(although we have abused notation a bit here by scaling to omit the constant $0.00926927$). Although the roots of this polynomial can be understood using work of Conrey, Farmer, and Imamoglu in \cite{ConreyFarmerImamoglu13}, this polynomial does not satisfy the Riemann hypothesis that $R_\Delta(X)$ does, so the work in \cite{OnoRolenSprung17} does not apply here.  Since $R_\Delta^-(X)$ still satisfies equations \eqref{fricke}, \eqref{Res1}, and \eqref{Res2}, Theorems \ref{eichshim1} and \ref{eichshim2} still apply to the zeta-polynomial
\begin{equation*} \begin{split}
Z_\Delta^-(s) = {\frac {{s}^{10}}{36288}}-{\frac {5\,{s}^{9}}{36288}}+{\frac {7\,{s}^{
8}}{2160}}-{\frac {367\,{s}^{7}}{30240}}+{\frac {833\,{s}^{6}}{8640}}-
{\frac {2137\,{s}^{5}}{8640}}\\
+{\frac {70841\,{s}^{4}}{90720}}-{\frac {
13193\,{s}^{3}}{11340}}+{\frac {403\,{s}^{2}}{360}}-{\frac {727\,s}{
1260}}.
\end{split} \end{equation*}
In particular, we have that
\[Z_\Delta^-(s) - Z_\Delta^-(1-s) =0.\]
Analogous statements hold for $Z_\Delta^+(s)$ as well.
\end{ex}

\begin{thm}\label{hilbertpoly}
Let $w\geq 2$ be even, let $R(X) \in \mathbb{C}[X]$ be any polynomial of degree at most $w$, and let $Z(s)$ be the polynomial satisfying \eqref{defofZ}. If $R(X)$ satisfies \eqref{fricke} and $Z(s)$ has integer coefficients with positive leading term, then $Z(s)$ is a Hilbert polynomial.
\end{thm}


This paper is organized as follows. In Section \ref{prelim}, we will review the relevant background related to period polynomials, zeta-polynomials, and Hilbert polynomials.  In Sections \ref{proofes}, \ref{proofes2}, and \ref{proofh}, we will prove Theorems \ref{eichshim1}, \ref{eichshim2}, and \ref{hilbertpoly}, respectively.


\section{Preliminaries}\label{prelim}

\subsection{Period polynomials of modular forms} \label{perpolybackground}
First we must define our notation and review the required background related to modular forms and their period polynomials; period polynomials give a context for Theorems \ref{eichshim1}, \ref{eichshim2}, and \ref{hilbertpoly} by providing natural applications of these results. For additional information, see, for example, the discussions in  \cite{KohnenZagier84} and \cite{ChoieParkZagier19}. 

Here we follow the standard notation: let $\mathbb{H}$ denote the upper half plane. For an even integer $k$ and $\gamma = \pm \begin{pmatrix}a&b\\ c&d\end{pmatrix}\in \PGL_2^+(\mathbb{R})$ we define the \emph{slash operator} $\mid_k$ for holomorphic functions $f:\mathbb{H}\rightarrow \mathbb{C}$ by
\[\left(f \mid_k \gamma\right)(\tau) := (ad-bc)^{k/2}(c\tau +d)^{-k}f\left(\frac{a\tau+b}{c\tau+d}\right).\]
If $N$ is a positive integer, we let $S_k(\Gamma_0(N))$ denote the space of cusp forms of weight $k$ on $\Gamma_0(N).$

Now we summarize the theory of period polynomials. Let $f\in S_k(\Gamma_0(N))$ be a cusp form of even weight $k$ and level $N$, and set $w:=k-2.$ The period polynomial associated to $f$ is given by
\[r_f(X) := \int_0^{\infty} f(\tau)(X-\tau)^{w}\,d\tau,\]
which is a polynomial in the space
\[\mathbb{V}_w := \{P \in \mathbb{C}[X] : \deg(P)\leq w\}.\]
We also define $r_f^+(X)$ and $r_f^-(X)$ to be the even and odd parts of $r_f(X)$, respectively, and note that $r_f^\pm (X)\in \mathbb{V}_w^\pm$ (where of course $\mathbb{V}_w^+$ and $\mathbb{V}_w^-$ are defined to be the set of even and odd polynomials of degree at most $w$, respectively).  There is an action of $\mathrm{PGL}_2^+(\mathbb{R})$ on $\mathbb{V}_w$ via the slash operator $\mid_{-w}.$

Note that if the cusp form $f$ is an eigenfunction of the Fricke involution $W_N = \begin{pmatrix}0&-1\\ N&0\end{pmatrix}$, i.e., $f\mid_k W_N = \varepsilon f$ for $\varepsilon \in \{\pm 1\}$, then it follows that $r_f$ (as well as $r_f^\pm$) satisfies
\begin{equation} \label{fricker}
r_f \mid_w \left(1 + \varepsilon W_N\right) = 0.
\end{equation}
(This fact can also be obtained using the functional equation of the $L$-function associated to $f$.)

On the other hand, if the modular form $f$ has level $N=1,$ then one can show that $r_f$ (as well as $r_f^\pm$) satisfies the Eichler-Shimura relations
\begin{align}
r_f|(1 + S) &= 0 \label{es1}\\
r_f|(1 + U+U^2) &= 0, \label{es2}
\end{align}
where
\[S := \begin{pmatrix}0&-1\\1&0\end{pmatrix} \qquad U := \begin{pmatrix} 1&-1\\ 1&0\end{pmatrix}.\]
Thus we define
\[\mathbb{W}_{w} := \{P \in \mathbb{V}_w : P\mid (1+S) = P\mid (1+U+U^2) = 0\}\] and note that $r_f\in \mathbb{W}_w.$ The following result of Eichler-Shimura illustrates the importance of the period polynomial $r_f(X).$
\begin{thm}
The map
\begin{align*} 
S_k(\SL_2(\mathbb{Z})) &\rightarrow \mathbb{W}_w^-\\
f &\mapsto r_f^-
\end{align*}
is an isomorphism. The map
\begin{align*}
S_k(\SL_2(\mathbb{Z})) &\rightarrow \mathbb{W}_w^+\\
f &\mapsto r_f^+
\end{align*}
is an injection whose image is a subspace of $\mathbb{W}_w^+$ of codimension 1.
\end{thm}

\subsection{Zeta-polynomials for modular form periods}
 
Let $f\in S_k(\Gamma_0(N)$ be a newform of even weight $k\geq 4.$ As discussed in the introduction, Ono, Rolen, and Sprung considered in \cite{OnoRolenSprung17} a reformulated version of the period polynomial
\[R_f(X) := (\sqrt{N}/i)^{k-1}r_f(X/i\sqrt{N}).\]
These polynomials serve as the inspiration for this work, so we note here that equation \eqref{fricker} gives
\[R_f(X) + \varepsilon i^wX^wR_f(1/X) = 0,\]
i.e., $R_f(X)$ satisfies equation \eqref{fricke}. Also, in the special case where $N=1$, the Eichler-Shimura relations \eqref{es1} and \eqref{es2} above  tell us that
\begin{align*} 
R_f(X) + (-iX)^w R_f(1/X)&= 0\\
R_f(X) + (-iX)^wR_f\left(\frac{X-i}{-iX}\right) + (-iX-1)^wR_f\left(\frac{-i}{-iX-1}\right) &= 0,
\end{align*}
i.e., $R_f(x)$ satisfies equations \eqref{Res1} and \eqref{Res2} when $N=1$ (with $\varepsilon=1$). As discussed in the introduction, the zeta-polynomials for modular form periods $Z_f(s)$ are given by
\[\frac{R_f(X)}{(1-X)^{w+1}} = \sum_{n\geq 0}Z_f(-n)X^n.\]

\subsection{Hilbert polynomials} Here we give the necessary background related to Hilbert polynomials; for more information, see \cite{Brenti98}. Fix a field $k$, let $\displaystyle{R = \bigoplus_{j \geq 0}R_j}$ be a graded $k$-algebra, and suppose that $R$ is standard (i,.e., that it can be finitely generated by elements of $R_1$). The \emph{Hilbert series} of $R$ is the formal power series \[\sum_{j\geq 0}\mathrm{dim}_k(R_j)X^j.\]
It is known that the Hilbert series can be written as
\[\frac{U(X)}{(1-X)^r} = \sum_{j\geq 0}\mathrm{dim}_k(R_j)X^j\]
for some positive integer $r$ and some polynomial $U(X),$ and it is also known that there exists a polynomial $P_R(X)\in \mathbb{Q}[X]$ such that \[P_R(j) = \mathrm{dim}_k(R_j)\] for all sufficiently large $j$.

Thus we make the following definition: a polynomial $H(X) \in \mathbb{Q}[X]$ is called a \emph{Hilbert polynomial} if there exists a standard graded $k$-algebra $R$ such that $H(X) = P_R(X).$ Work of Brenti \cite{Brenti98} investigates which polynomials are Hilbert polynomials, and how to measure ``how far'' a polynomial is from being Hilbert.  Along the way, Brenti proves the following useful results.

\begin{thm}[Theorems 3.5 and 3.14 of \cite{Brenti98}] \label{Bthm}
Let $H(X)\in \mathbb{Z}[X]$ be a polynomial with positive leading term.
\begin{itemize}
\item There exists $M \in \mathbb{N}$ such that $H(X+j)$ is a Hilbert polynomial for any $j\geq M.$
\item If $H(X)$ is a Hilbert polynomial then $H(X+1)$ is a Hilbert polynomial.
\end{itemize}
\end{thm}


\section{Proof of Theorem \ref{eichshim1}} \label{proofes}

Let $w\geq 2$ be even, let $R(X) = \sum_{j=0}^wa_jX^j \in \mathbb{C}[X],$ and let $Z(s)$ be the polynomial satisfying
\[\frac{R(X)}{(1-X)^{w+1}} = \sum_{n\geq 0}Z(-n)X^n.\]
In order to better understand the relationship between $R(X)$ and $Z(s)$, we use Newton's Binomial Theorem, which says that
\[\frac{1}{(1-X)^{w+1}} = \sum_{n\geq 0}\binom{w+n}{n}X^n.\]
Thus we have
\begin{align*}
\sum_{n\geq 0}Z(-n)X^n &= \frac{R(X)}{(1-X)^{w+1}}\\
&=  \left( \sum_{j=0}^wa_jX^j \right)\left( \sum_{n\geq 0}\binom{w+n}{n}X^n \right)\\
&= \sum_{n\geq 0} \left(\sum_{j=0}^{w}a_j\binom{w+n-j}{w}\right)X^n
\end{align*}
so we can now express $Z(s)$ explicitly by
\[Z(s) = \sum_{j=0}^{w}a_j\binom{w-s-j}{w}.\]

Now, to prove Theorem 1, we suppose equation \eqref{fricke}, i.e., that 
\[a_j + \varepsilon i^w a_{w-j}=0\]
for $0\leq j\leq w.$ Thus
\begin{align*}
Z(1-s) &= \sum_{j=0}^{w} a_j\binom{w-1+s-j}{w} = \sum_{j=0}^{w} a_{w-j}\binom{w-1+s-(w-j)}{w}\\
&= \sum_{j=0}^{w} -\varepsilon i^w a_{j}\binom{-1+s+j}{w}  = \sum_{j=0}^{w} -\varepsilon i^w a_{j}\binom{w-s-j}{w}\\
&= -\varepsilon i^w Z(s),
\end{align*}
as desired.

\qed


\section{Proof of Theorem \ref{eichshim2}} \label{proofes2}

First, note that the first statement of Theorem \ref{eichshim2} follows from Theorem \ref{eichshim1} by letting $\varepsilon=1$. To complete the proof of Theorem \ref{eichshim2}, we note that equation \eqref{Res2} tells us that for any positive integer $n$, we have that
\begin{equation} \label{RIes2}
\begin{split}
\frac{1}{2\pi i}\int_\gamma \frac{R(z)}{(1-z)^{w+1}}z^{-(n+1)} \,dz + \frac{1}{2\pi i}\int_\gamma \frac{(-iz)^wR\left(\frac{z-i}{-iz}\right)}{(1-z)^{w+1}}z^{-(n+1)} \,dz\\ + \frac{1}{2\pi i}\int_\gamma \frac{(-iz-1)^wR\left(\frac{-i}{-iz-1}\right)}{(1-z)^{w+1}}z^{-(n+1)} \,dz = 0,
\end{split}
\end{equation}
where $\gamma$ is a small circle with center 0 (oriented counter-clockwise). Our proof will follow by interpreting each integral of equation \eqref{RIes2}; note that by Cauchy's integral formula, the first integral is
\[\frac{1}{2\pi i}\int_\gamma \frac{R(z)}{(1-z)^{w+1}}z^{-(n+1)} \,dz = Z(-n).\]

Next, we note that the second integral is (by applying equation \eqref{Res1} and then setting $x = -iz/(z-i)$)
\begin{align*}
\frac{1}{2\pi i}\int_\gamma \frac{(-iz)^wR\left(\frac{z-i}{-iz}\right)}{(1-z)^{w+1}}z^{-(n+1)} \,dz &= \frac{-1}{2\pi i}\int_\gamma \frac{(-iz)^w\left(\frac{z-i}{z}\right)^w R\left(\frac{-iz}{z-i}\right)}{(1-z)^{w+1}}z^{-(n+1)} \,dz\\
&= \frac{-(-i)^w}{2\pi i}\int_\gamma \frac{R\left(\frac{-iz}{z-i}\right)(z-i)^{w+2} }{(1-z)^{w+1}}z^{-(n+1)} \,\frac{dz}{(z-i)^2}\\
&= \frac{(-i)^w}{2\pi i}\int_{\gamma_0} \frac{R(x)}{(1-x)^{w+1}} \frac{(1-x)^{w+1} (x+i)^{w+1}(x+i)^{n+1}}{(x+i-ix)^{w+1}(x+i)^{w+2}(ix)^{n+1}} \,dx\\
&= \frac{(-i)^w}{2\pi i}\int_{\gamma_0} \frac{R(x)}{(1-x)^{w+1}} \frac{(1-x)^{w+1} (x+i)^{n}}{(x+i-ix)^{w+1}(ix)^{n+1}} \,dx\\
&= (-i)^w\sum_{m=1}^{n+1} a_{-m}Z(1-m)
\end{align*}
where $\gamma_0$ is a small circle with center 0 and
\[\frac{(1-x)^{w+1} (x+i)^{n}}{(x+i-ix)^{w+1}(ix)^{n+1}} = \sum_{m=-n-1}^\infty a_mx^m.\]

Finally, by Cauchy's integral formula and equation \eqref{Res1}, the third integral is the coefficient of $z^n$ in the expansion for
\[\frac{(-iz-1)^wR\left(\frac{-i}{-iz-1}\right)}{(1-z)^{w+1}} = \frac{-R(z-i)}{(1-z)^{w+1}}\]
about 0. Thus we define a function $g$ by
\[g(z-i) = \frac{-R(z-i)}{(1-z)^{w+1}} = \frac{-R(z-i)}{(1-z+i)^{w+1}}\cdot \frac{(1-z+i)^{w+1}}{(1-z)^{w+1}}.\]
Then
\begin{align*}
g(z) &= \frac{-R(z)}{(1-z)^{w+1}}\cdot \frac{(1-z)^{w+1}}{(1-i-z)^{w+1}}\\
&= \left(-\sum_{n\geq 0}Z(-n)z^n\right)(1-z)^{w+1}\left(\sum_{m \geq 0}\binom{m+w}{w}\frac{z^m}{(1-i)^{m+w+1}}\right)\\
&= \sum_{k\geq 0} \left[ \sum_{m=0}^k  \binom{m+w}{w}\frac{1}{(1-i)^{m+w+1}} \sum_{j=0}^{k-m} \binom{w+1}{j}(-1)^{j+1} Z(m+j-k) \right]z^k.\\
\end{align*}
Now, set $b_k$ to be the expression inside the brackets above, so that $g(z) = \sum_{k\geq 0}b_kz^k.$ Since the radius of convergence of this series is $\sqrt{2},$ we may substitute to find
\begin{align*}
g(z-i) &= \sum_{k\geq 0}b_k(z-i)^k = \sum_{n\geq 0}\sum_{k\geq 0} b_{k+n}\binom{k+n}{n}(-i)^kz^n\\
&= \sum_{n\geq 0} \left[ \sum_{k\geq 0}\sum_{m=0}^{k+n} \sum_{j=0}^{k+n-m} \binom{k+n}{n} \binom{m+w}{w} \binom{w+1}{j} \frac{(-1)^{j+1}(-i)^k }{(1-i)^{m+w+1}} Z(m+j-k-n) \right]z^n.
\end{align*}
Thus the third integral is $\sum_{k\geq 0}\sum_{m=0}^{k+n} \sum_{j=0}^{k+n-m} \binom{k+n}{n} \binom{m+w}{w} \binom{w+1}{j} \frac{(-1)^{j+1}(-i)^k }{(1-i)^{m+w+1}} Z(m+j-k-n),$ completing the proof. \qed


\section{Proof of Theorem \ref{hilbertpoly}} \label{proofh}

Suppose for the sake of contradiction that $Z(s)$ is not a Hilbert polynomial. By the first part of Theorem \ref{Bthm}, there exists some $M\in \mathbb{N}$ such that $Z(s+j)$ is Hilbert for any $j\geq M$ (and we may suppose without loss of generality that $M$ is minimal, i.e., that $Z(s+M)$ is Hilbert and $Z(s+M-1)$ is not Hilbert).

Note that by Theorem \ref{eichshim1} \[Z(s+M) = -\varepsilon i^w Z(1-M-s)\] is Hilbert, so the second part of Theorem \ref{Bthm} tells us that
\[-\varepsilon i^wZ(2-M-s) = Z(M+s-1)\]
is also Hilbert.  This is a contradiction. \qed

\bibliography{/Users/jameson/Dropbox/LaTeX/master}
\bibliographystyle{alpha}

\end{document}